\documentclass{conm-p-l}

\usepackage{amssymb}
\usepackage{epsfig}
\usepackage{amsfonts}
\usepackage{amsmath}
\usepackage{amssymb}
\newcommand{\intersect}{\cap}
 \newcommand\supp{\mathop{\rm supp}} 
\newtheorem{theorem}{Theorem}[section]

\newcommand*{\pd}[2]{\frac{\partial #1}{\partial #2}}

\begin{document}
%amsmath
%\begin{frontmatter}




%\title{On the control of the interface for the Muskat problem}
\title[Exact solutions to a Muskat problem]{Exact solutions to a Muskat problem with line distributions of sinks and sources}
\author{L.~Akinyemi}
%\author[Tanya1,Tanya2,Tanya3]{T.V. Savina\thanksref{bby}\corauthref{cor}}
\author{T.V.~Savina}
\address{Department of Mathematics,  Ohio University, Athens, OH 45701, USA}
\address{Condensed Matter and Surface Science Program, Ohio University, Athens, OH 45701, USA}
\address{Nanoscale \& Quantum Phenomena Institute, Ohio University, Athens, OH 45701, USA. }
\email{savin@ohio.edu}
%\thanks[bby]{Research of this author was supported in part by OU Research Challenge Program, award \# RC-09043.} and
%\ead{savin@ohio.edu} \ead[url]{http://www.ohio.edu/people/savin/}
%\corauth[cor]{Corresponding author.}
%\author[Alex1,Alex2]{A.A.~Nepomnyashchy\thanksref{sa}}
\author{A.A.~Nepomnyashchy}
%\ead{nepom@techunix.technion.ac.il} \ead[url]{http://www.math.technion.ac.il/Site/people/process.php?id=710}
%\thanks[sa]{Research of this author was supported in part  by the European Union via FP7 Marie Curie scheme Grant 
%PITN-GA-2008-214919 (MULTIFLOW).} 
\address{Department of Mathematics, Technion - Israel Institute of Technology,
 Technion City,  Haifa,  32000, Israel}
\address{Minerva Center for Nonlinear Physics of Complex Systems, Technion - Israel Institute of Technology,
 Technion City,  Haifa,  32000, Israel}

\begin{abstract} 
The evolution of a two-phase Hele-Shaw problem, a Muskat problem, under assumption of a negligible surface tension is 
 considered. 
We use the Schwarz function approach and allow the sinks and sources to be  line distributions with  disjoint supports located in  the exterior and the interior domains, a two-phase mother body. 
We give examples of exact solutions when the interface belongs to  a certain family of algebraic curves, defined by the initial shape of the boundary, and the cusp formation does not occur. 
\end{abstract}
%{\bf Key words:} 

%\begin{keyword}
%\end{keyword}

\begin{date}
\date{}
\end{date}

% END topmatter
%\end{frontmatter}
\maketitle

\section{Introduction}

A Muskat problem \cite{mus} describes an evolution of an interface between two immiscible fluids, `oil' and `water',  in a Hele-Shaw cell or in a porous medium.
 Hele-Shaw free boundary problems have been extensively studied over the last century (see \cite{Vas2009}, \cite{Vas2015} and references therein). There are two classical formulations of the Hele-Shaw problems: the one-phase problem, when one of the fluids is assumed to be viscous while the other is effectively inviscid (the pressure there is constant), and the two-phase (or Muskat) problem. 
For the latter problem much less progress has been made than for the former one. Concerning the two-phase problem, we should mention works \cite{howison2000}-\cite{crowdy2006}. Specifically,
Howison \cite{howison2000} has obtained several simple solutions including the traveling-wave solutions and the stagnation point flow. In \cite{howison2000}, an idea of a method for solving some two-phase problems was  proposed and used to reappraise the Jacquard-S\'eguier solution \cite{JS}.
Global existence of solutions to some specific two-phase problems was considered in 
\cite{FT}-\cite{YT}.
Crowdy \cite{crowdy2006} presented an exact solution to the Muskat problem for the elliptical initial interface between two fluids of different viscosity. In \cite{crowdy2006} it was shown that an  elliptical inclusion of one fluid remains elliptical when placed in a linear ambient flow of another fluid.
 This paper concerns exact solutions to the Muskat problem, extending the results obtained in
\cite{crowdy2006} to other types of inclusions.  
The main difficulty of the two-phase problems is the fact that the pressure on the interface is unknown. However, if we assume that the free boundary remains within the family of curves, specified by the initial shape of the interface separating the fluids (which is feasible if the surface tension is negligible) the problem is drastically simplified. This study is devoted to the situations
when the evolution of the interface is controlled by a special choice of  sinks and sources.
The suggested method allows to obtain exact solutions for a certain class of curves for which the Schwarz function can be computed.

The mathematical formulation of the problem is as follows.
Let $\Omega _2 (t) \subset {\mathbb R}^2$ with a boundary $\Gamma (t)$ at time $t$ be a 
simply-connected bounded domain occupied by a fluid with a constant viscosity $\nu _2$,   and let $\Omega _1 (t)$ be the region  ${\mathbb R}^2\setminus {\bar \Omega}_2(t)$ occupied by a different fluid of viscosity $\nu _1$.
Consider the two-phase Hele-Shaw problem forced by sinks and sources:
\begin{equation}\label{1}
{\bf v}_j=-k_j\nabla p_j, \qquad j=1,2, 
%\qquad \Delta p=0 \quad\mbox{in}\quad \Omega (t).
\end{equation}
where the pressure $p_j$ is a harmonic function almost everywhere in the region $\Omega _j (t)$, satisfying boundary conditions
\begin{eqnarray}\label{2}
p_1(x,y,t)=p_2(x,y,t) \quad \mbox{on} \quad \Gamma (t),\\
-k_1\pd{p_1}{n}=-k_2\pd{p_2}{n} =v_n \quad \mbox{on} \quad \Gamma(t).\label{3}
\end{eqnarray}
Here ${\bf v}_j$ is a velocity vector of fluid $j$, $k_j=h^2/12\nu _j$, and $h$ is the gap width of the Hele-Shaw cell.
Equation (\ref{2}) states the continuity of the pressure under the assumption of negligible
 surface tension.  
 Equation (\ref{3}) means that the normal velocity of the boundary itself coincides with the normal velocity of the fluid
 at the boundary.

The free boundary $\Gamma (t)$ moves due to the sources and sinks  
 located in both regions.
We adopt a natural physical assumption that the fluid flux generated by the system of sources and sinks is finite. That allows no more than the logarithmic growth of the fluid pressure near a point source/sink or at  infinity,
% and the pressure at any point of the plane $(x_a,y_a)$   has  no more than a logarithmic singularity, 
%p(x,y,t)=-\frac{Q}{2\pi}\log \sqrt{x^2+y^2}\quad\mbox{as}\quad \sqrt{x^2+y^2}\to\infty . \label{4}
\begin{equation}
|p_j(x,y,t)|\le |-\frac{Q_a(t)}{2\pi k_j}\log \sqrt{(x-x_a)^2+(y-y_a)^2}\,|
%\quad\mbox{as}\quad \sqrt{x^2+y^2}\to\infty 
, \label{4m}
\end{equation}
where $Q_a(t)$ is the strength of the source/sink. We would like to stress that this physical assumption is not a restriction of the suggested method. In Section \ref{sec:ex}, for the sake of curiosity, we obtain a solution, whose far field flow is linear. The latter type of flow was  obtained by  Crowdy \cite{crowdy2006}.

In this study, we allow  the supports of sinks and sources to not only be points, but lines/curves as well, which could essentially change the dynamics of the evolution of the interface. A similar approach was used for the one-phase problems  \cite{dmb}, \cite{external}. 
In the case of the interior problem, this approach was     
motivated by the fact that during extraction through a point sink located within a viscous fluid, the free boundary  is unstable, and the solution breaks down before all the fluid is extracted due to the formation of cusps, except for the situation of a circular boundary with a sink in the center.  In \cite{dmb} it was shown that a choice of  sinks with line distributions linked to the initial shape occupied by the viscous fluid allows  to enlarge the class of domains, from which the viscous fluid can be completely extracted without a cusp formation.
Analogously,
for the exterior problem with uniform extraction at infinity 
Howison \cite{HowBubble} has proven that the elliptical bubbles are the only finite bubbles which exist for all times and whose boundary crosses all points initially outside the bubble.
 In all other cases, the solution either fails to exist in a finite time or the solution has some points on the interface that have a finite limit as time approaches infinity, so some fluid is ``left behind'' \cite{di}.  
In the recent work \cite{external}, it was shown that if a point sink at infinity is replaced with a specific line distribution of sinks in the exterior region, then the evolution changes and 
it is possible to find other than elliptical shapes for which in the course of growing, the boundary of the air bubble  crosses all points outside the bubble.

In the present paper we introduce the notion of the two-phase mother body, which generalizes the
dynamical (one-phase) mother body used  in \cite{dmb} and \cite{external}.
% the notion of dynamical (one-phase) mother body used in \cite{dmb} and \cite{external} is developed for the case of  the two-phase Hele-Shaw problem, which 

We  would like to acknowledge 
the results obtained by Karp related to the unbounded quadrature domains, including the asymptotic behavior of the boundary in $\mathbb{R}^2$ \cite{karp1} and the  connections between the generalized Newtonian potential and the unbounded quadrature domains in  $\mathbb{R}^n$, $n\ge 3$ \cite{karp2}.

It is also worth  mentioning
the recent development in the two-phase quadrature domain theory \cite{shah}, \cite{gardiner}.
In the spirit of the latter theory the problem in question could be reformulated as follows.
Let $u(x,y)$ be a continuous across $\Gamma (t)$ function, such that  $u(x,y)\chi _{[{\Omega} _1]}=p_1$ and $u(x,y)\chi _{[\Omega _2]}=p_2$, where
$$
\chi _{[\Omega _j]}=\left \{
\begin{aligned}
1, \quad \mbox{if} \quad (x,y)\in\Omega _j\cup\Gamma\\
0, \quad\mbox{if}  \quad (x,y)\not\in \Omega_j \cup\Gamma.\\
\end{aligned}
\right .
$$
Find a solution to the problem
\begin{eqnarray}\label{Q1}
\Delta u= \mu_1 +\mu_2 \quad \mbox{in} \quad \mathbb{R}^2  ,\\
%u=0 \quad \mbox{on} \quad \Gamma (t), \label{Q2}\\ 
-k_j\pd{u}{n}\chi _{[\Omega _j]}\to v_n \quad \mbox{as} \quad (x,y)\to \Gamma (t),\label{Q3}
\end{eqnarray}
where  $\mu_j(t)$ are time dependent distributions with $\supp\, \mu _j(t)\subset \Omega _j$, $\supp\,\mu_1\intersect \supp\,\mu _2=\varnothing$.

The structure of the paper is as follows.    In Section \ref{prelim} we reformulate the problem in terms of the Schwarz function and indicate the steps of the suggested method of finding exact solutions. Section \ref{sec:ac} is devoted to the two-phase mother body in the context of the Muskat problem.
Examples of the exact solutions are given in Section \ref{sec:ex}, and conclusions are given in Section \ref{sec:concl}.

\section{The Schwarz function method of finding exact solutions for the Muskat problem}\label{prelim}

As mentioned above, the 
evolution of the interface separating the fluids is determined by the distributions of sinks and sources, which in the absence of the surface tension, could be chosen in such a way that keeps 
$\Gamma(t)$ within a  family of curves defined by $\Gamma (0)$. 
For what follows, it is convenient to reformulate problem   (\ref{1})--(\ref{4m}) in terms of the Schwarz function $S(z,t)$ of the curve $\Gamma (t)$ \cite{davis}--\cite{shapiro}.
%, \cite{laguna},  \cite{nonloc}, 
% \cite{laguna}, \cite{bihar}, 
This function for
 a real-analytic curve  $\Gamma :=\{ g(x,\,y,\,t)=0\}$  is defined as a solution 
to the equation
%\begin{equation}\label{E:999}
$g\left ((z+\bar z)/2, \, (z-\bar z)/2i,\, t \right ) =0$
%\end{equation}
with respect to $\bar z$.
%Here $g(x,\,y,\,t)$ is a given function that can be  extended into a holomorphic function, for example, a polynomial.
Such (regular) solution exists in some neighborhood $U_{\Gamma}$ of the curve $\Gamma$,
if the assumptions  of the implicit function theorem are satisfied \cite{davis}.
Note that if $g$ is a polynomial, then the Schwarz function is continuable into $\Omega _j$,
generally as a multiple-valued analytic function with a finite number of algebraic singularities (and poles).
In $U_{\Gamma}$, the normal velocity, $v_n$, of  $\Gamma (t)$ can be written in terms of the Schwarz function \cite{howison92},
%\begin{equation}\label{velocity}
$v_n=-i\dot S(z,t)/\sqrt{4\partial _z S(z,t)}$.
%\end{equation}
%To complete reformulation of the problem, we introduce the complex potentials $W_j=p_j-i\psi _j$, which are analytic multiple-valued functions 
%defined in $\Omega _j (t)$ in the neighborhood of $\Gamma (t)$.  Here $\psi _j$ is a stream function.
% Remark that, 
%$v_x^j=\pd{\psi _j}{y}$ and $v_y^j=-\pd{\psi _j}{x}$,
%where ${\bf v}^j=(v_x^j,v_y^j)$ is a velocity field of an incompressible potential flow near $\Gamma (t)$. Thus, 
%$\pd{W_j}{z}=\pd{p_j}{x}-i\pd{p_j}{y}=-v_x^j+iv_y^j.$

Let $\tau$ be an arclength along $\Gamma (t)$, $\psi _j$ be a stream function, and $W_j=p_j-i\psi _j$ be the complex potential defined  on $\Gamma (t)$ and in $\Omega _j (t)\cap U_{\Gamma}$, $j=1,\,2$. Following \cite{cummings}-%
%\cite{KhMP}
\cite{lacey},  taking into account the Cauchy-Riemann conditions in the $(n,\tau )$ coordinates, for the derivative of $W_j(z,t)$ with respect to $z$ on $\Gamma (t)$ we have
\begin{equation}\label{main1}
\partial _z {W_j}=\frac{\partial _{\tau}W_j}{\partial _{\tau}z}=
%\frac{\partial _{\tau}(p_j-i\psi _j)}{\partial _{\tau}z}=
\frac{\partial _{\tau}p_j+i\partial _{n}p_j}{\partial _{\tau}z}=
\frac{\partial _{\tau}p_j-iv_n/k_j}{\partial _{\tau}z}.
\end{equation}
Expressing $\partial _{\tau}z$ in terms of the Schwarz function, $\partial _{\tau}z=(\partial _z S(z,t))^{-1/2}$, we obtain 
\begin{equation}\label{main1m}
\partial _z {W_j}=\partial _{\tau}p_j\sqrt{\partial _z S}
-\frac{\dot {S}}{2k_j}.
\end{equation}
Here $\partial _z {W}_j\equiv \pd{W_j}{z}$,  $\partial _{\tau} {z}\equiv \pd{z}{\tau}$ etc.
%Note that according to  (\ref{2}) real part of $W(z,t)$ is constant on $\Gamma (t)$ and the latter equation, 
Since $p_1=p_2$ on $\Gamma (t)$,  equation \eqref{main1m} implies
\begin{equation}\label{main}
\partial _z {W_1}+\frac{\dot {S}}{2k_1}=\partial _z {W_2}+\frac{\dot {S}}{2k_2}
=\partial _{\tau}p_j\sqrt{\partial _z S}.
\end{equation}
To keep $\Gamma (t)$ in a certain family  of curves defined by $\Gamma (0)$, for example, in a family of ellipses, we  assume that $p_j$ on $\Gamma (t)$ is a function of time only. In that case the problem is simplified drastically, and on $\Gamma (t)$ we have
\begin{equation}\label{main2m}
\partial _z {W_j}=
-\frac{\dot {S}}{2k_j}\qquad j=1,2.
\end{equation}
Taking into account that $\dot S$ can be continued off of $\Gamma (t)$,
%is defined in  $U_{\Gamma}$, 
each equation \eqref{main2m} can be continued off of $\Gamma$ into the corresponding
$\Omega _j$, where $W_j$ is a multiple-valued analytic function. 

Note that equations \eqref{main2m} indicate that the singularities of $W_1$, $W_2$ and the Schwarz function are linked. Therefore, as we show below, those singularities in some cases can be used to  control the interface between the fluids. Thus, the problem reduces to finding the distributions $\mu _1 (t)$ and $\mu_2 (t)$, that keep $\Gamma (t)$ in a family of curves generated by $\Gamma (0)$.  
The latter problem can be viewed as a generalization of a classical problem of electrostatics: find a two-charge system that yields the (desired) zero potential on a conducting plate.

To find the exact solutions,
suppose that at $t=0$ the interface is an algebraic curve, $\sum _{k=0}^n a_k(0) x^{k-n}y^n=0$, with the
Schwarz function $S(z,a_k^0)$. Assume  that during the course of evolution the Schwarz function of the interface $S(z,a_k(t))\equiv S(z,t)$ is such that
$S(z,a_k(0))=S(z,a_k^0)$. 

The steps of the method are
\newline
1) Compute $\dot S (z, t)$, locate its singularities, and define their type.
\newline
2) Using equations \eqref{main2m} find  preliminary expressions for $\partial _z W_j$ and,
  by putting restrictions on the coefficients $a_k(t)$,
eliminate their terms that involve non-integrable singularities.
\newline
3) Find the quantities $W_j$ by integrating \eqref{main2m} with respect to $z$. 
\newline
4) Compute the quantities $p_j$ by taking the real parts of $W_j$.
\newline
5) Evaluate the quantities $p_j$ on the interface to determine the independent of $z$ function of integration from the step 3).
\newline
6) Compute the two-phase mother body.
\newline
We comment that the steps 1)-5) are straight-forward, and the step 6) is discussed below.

\section{A two-phase mother body}\label{sec:ac}

Generally, the complex potentials $W_j$ are multiple-valued functions in $\Omega _j$.
To choose a branch for each of  these functions, one has to introduce the cuts, that serve as 
supports for the distributions of sinks and sources. The union of these distributions  $\mu _1 (t)$, $\mu _2(t)$ with  disjoint supports (see formula \eqref{Q1}) and integrable densities, which allows a smooth evolution of the interface, is called below
{\it a two-phase mother body}. 
The notion of  a mother body comes from the potential theory
\cite{shah}, \cite{gardiner}, 
%\cite{Zid1}-
%, \cite{Zid2},  \cite{GuSa2},   
%\cite{ts},  
\cite{Gu1}-\cite{sss2005}.  
 The supports of these distributions consist of   sets of arcs  and/or points and  do not bound any two-dimensional subdomains in 
$\Omega _j(t)$, $j=1,2$. 
Each cut included in the support of  $\mu _j(t)$ is contained in the domain $\Omega _j (t)$, and the limiting values of
the pressure on each side of the cut are equal. 
The value of the source/sink
density on the cut is equal to the jump of the normal derivative $\partial _n p_j$ of the pressure $p_j$.
To ensure that the total flux through the sources/sinks is finite,
all of the singularities of the function $W_j$ must have no
more than logarithmic growth.

If $\Gamma (t)$ is an algebraic curve, then the singularities of  $W_j$ are either poles or algebraic singularities. Thus,
 each  cut  originates from 
an algebraic singularity  $z_a(t)$ of the potential $W_j$.
To ensure that  the limiting values of $p_j$ on both sides of each cut are
%$p={\rm Re}\, W$ 
equal, using the terminology of algebraic topology (\cite{algtop}, p. 21), we proceed as follows. We fix a point $z_b\in\Gamma (t)$, the base point, and consider a fundamental group of loops $l\subset \Omega _j (t)$, having $z_b$ as their starting and the terminal point, and surrounding the singular point $z_a(t)$, the group $\pi _1(\Omega _j \setminus$sing $(W_j), z _a)$. 
A multiple-valued function  $W_j$ varies along $l$. We denote its variation by  ${\rm var}_l\,W_j$, and the real part of its variation by  ${\rm var}_l\,p_j$.
 %as a real part of this variation,
%${\rm var}_l\,W$,  
%${\rm var}_l\,p_j=\Re \{{\rm var}_l\,W_j\}$. 
Then, the zero level sets 
%  ${\rm var}_l\,p_j$ along any loop $l$,
%\subset \Omega (t)$ with the base point (the beginning and the ending point of the loop) $z_p\in\Gamma (t)$ surrounding the singular point $z_0(t)$.
%Thus, our goal is: for each $z_0(t)$  find a line or a curve having $z_0(t)$ as its endpoint and satisfying the equation  
${\rm var}_l\,p_j=0$ describe the location of the desired cuts.

%Obviously, all lines ${\rm var}_l\,p=0$ could be described as integral curves of some vector field.  To find this field let's fix time $t$. For each $t$ the boundary of $\Omega (t)$ is a zero level set of $p$. Thus, the desired  vector field at each point in $\Omega (t)$  must be orthogonal to the level curves  $p(x,y)=const$, and therefore, its direction at each point must coincide with  $\nabla p (x,y)$. Hence the lines  ${\rm var}_l\,p=0$ coincide with the lines $\Im [W]=const$, which are the stream lines. For the construction of an external mother body we need to take the lines that start at $z_0(t)$.

%Below we study the local structure of
%the cuts in a  neighborhood of each singularity $z_0(t)$. In the general case this analysis gives us only possible directions of each cut, but in some symmetric cases it gives a sufficient information  to describe the global structure of the cuts.
% in a similar fashion as in \cite{dmb}. 

Typically, the algebraic singularities are not stationary, that is $\dot z_a\ne 0$.
The location of $z_a(t)$ is determined by  the Schwarz function.
The theorem below states the uniqueness of the direction of the cut at a non-stationary singularity $z_a(t)$ in general position. 
%Equation (\ref{main}) suggests that 
% $S\left( z, t\right) $  in $\Omega _j (t)$. 
%We shall carry
%out this investigation under the assumption of the generic position. In other words, we
The latter means that the singularity $z_a(t)$  appears  from a finite
regular characteristic point of the 
%manifold $\Gamma _{\mathbb{C}}$, 
complexification of the boundary $\Gamma (t)$, and 
the tangency between this singular point and the corresponding
characteristic ray is quadratic. Under such requirements the function $%
S\left( z,t\right) $  at   $z_a(t)$ has  the square root
type singularity:
\begin{equation}\label{eq6}
S\left( z,t\right) =\Phi\left( z,t\right)\,\sqrt{z-z_a(t)} +\Psi\left(
z,t\right).
\end{equation}
Here $\Phi\left( z,t\right) $ and $\Psi\left( z,t\right) $ are regular functions of $z$ in a
neighborhood of the point $z_a(t)$, and $\Phi\left( z_a(t),t\right)\ne 0$.
The following theorem describes 
%possible initial conditions for the sought integral curves.
restrictions on the branch cuts in terms of their admissible slopes in the neighborhood of $z_a(t)$.
\begin{theorem}\label{cuts}
Let $z_a$ be a singular point of the complex potential $W_j$ located in $\Omega _j (t)$, $j=1,2$, such that $\dot z_a\ne 0$.
Then, under the assumption of general position \eqref{eq6},
%\newline 
% that is, $z_0$ is a moving singular point, 
the direction of the cut, on which ${\rm var}_l\,p_j=0$ near this point, is  uniquely defined by the formula
\begin{equation}\label{dir}
\varphi=\pi-2(\arg [\Phi\left( z_a(t),t\right)]+ \arg [\dot z_a]) +2\pi k,\quad k=0,\pm1, \pm 2 ... .
\end{equation} 
\end{theorem}
%$\cos \left( \frac{\varphi }2+\theta \right) =0$;
% where $\theta _0$ is an argument of $\Phi\left( z_0(t),t\right) $ and $$; 

{\bf Proof.}
We start with representation (\ref{eq6}) dropping the regular part, $\Psi\left( z,t\right) $, in it and
expanding  the function
$\Phi\left( z,t\right) $  into the Taylor series with respect to $z$ at the
point $(z_a(t),t)$, 
%say $\Phi(z,t)=\sum\limits _{j=0}^{\infty}c_j(t)\left( z-z_0(t)\right) ^j$. 
\begin{equation}
\label{eq67}S\left( z,t\right) =\sqrt{z-z_a(t)}\,\,\sum\limits _{m=0}^{\infty}c_m(t)\left( z-z_a(t)\right) ^m.
\end{equation}
The time derivative of the Schwarz function  (\ref{eq67}) has the form:
%$$
%\partial _t S\left( z,t\right) =-\frac{\partial _t z_0(t)}{2  \sqrt{z-z_0}}\,\sum\limits_{j=0}^%
%\infty c_j (2j+1)\left( z-z_0\right) ^j .
%$$
%\begin{equation}
%\begin{aligned}
%& \dot S\left( z,t\right) =-\frac{\dot z_0}{2  \sqrt{z-z_0}}\,
%\sum\limits_{j=0}^{\infty} c_j \left( z-z_0\right) ^j+ 
% \sqrt{z-z_0}\,\sum\limits_{j=0}^{\infty} \dot c_j \left( z-z_0\right) ^j\\
%& +\dot z_0 \sqrt{z-z_0}\,\sum\limits_{j=0}^{\infty} \Bigl [\pd{c_j}{z_0}-(j+1)c_{j+1}\Bigr ] \left( z-z_0\right) ^j,
%\end{aligned}
%\end{equation}
%where dot denotes the partial derivative with respect to $t$, while $\partial _t$ is the total derivative.
\begin{equation}
 \dot S \left( z,t\right) =-\frac{\dot z_a}{2}
\sum\limits_{m=0}^{\infty} c_m \left( z-z_a\right) ^{m-1/2}+ 
 \sum\limits_{m=0}^{\infty} \Bigl ( \dot c_m-\dot z_a (m+1) c_{m+1}\Bigr ) \left( z-z_a\right) ^{m+1/2},
\end{equation}
therefore, formula (\ref{main2m}) implies:
$$
W_j \left( z,t\right) =\dot z_a c_0 (z-z_a)^{1/2}\,\Bigl (\frac{1}{2k_j}+\xi_1(z,t) \Bigr )-\dot c_a  (z-z_a)^{3/2}\,\Bigl (\frac{1}{3k_j}+\xi_2(z,t)\Bigr ),
$$
where $\xi _1 $ and $\xi _2$ are regular functions near $z_a$ vanishing at
this point,
and $c_0=\Phi\left( z_a(t),t\right)$. 
Hence,
%$$
%p = \Re \, \Bigl \{ \dot z_0c_0 \sqrt{z-z_0}\,(\frac{1}{2} +s_1(z,t))-\dot c_0  (z-z_0)^{3/2}\,(\frac{1}{3}+s_2(z,t)) \Bigr \},
%$$
 the variation of the pressure along the loop $l$ is
\begin{equation}\label{strW}
{\rm var}_l\,p_j =\Re \,\Bigl \{ \dot z_a c_0 (z-z_a)^{1/2}\Bigl (\frac{1}{k_j}+2\xi_1(z,t)\Bigr )-\dot c_0  (z-z_a)^{3/2}\,\Bigl (\frac{2}{3k_j}+2\xi_2(z,t)\Bigr ) \Bigr \}.
\end{equation}
Consider a small neighborhood of the point $z_a$, and set there    $z=z_a+\rho e^{i\varphi }$, assuming that $\rho$ is small.
Since $\dot z_a \ne 0$, the second term in (\ref{strW}) is small with respect to the first term and therefore should be dropped. 
Setting the principal part of ${\rm var}_l\,p$ to  zero, 
we have
\begin{equation}\label{eq8}
{\rm var}_l\,p_j =\frac{1}{k_j}|\dot z_a||c_0|\sqrt{\rho }\, \Re \,
\,\Bigl \{ \exp i\left ( \frac{\varphi }{2}+\pi n+\arg \dot z_a+\arg c_0\right) \Bigr \}
%\left( 1+2\xi _1 \left( z_0+\rho e^{i\varphi }\right) \right) \Bigr \} 
=0,
\end{equation}
where $n=0,\pm 1,\pm 2, \dots$.
%$R_0$ and $\theta _0$ are the modulus and phase of $c_0$, $c_0=R_0 e^{i\theta _0 }$.

Equation \eqref{eq8} implies formula \eqref{dir}, which finishes the proof.

%Let us consider this equation as an equation with respect to $\varphi $ for
%small values of $\rho $. Then, 
%Since the function $s_1 \left( z_0+\rho
%e^{i\varphi }\right) $ is of order $O\left( \rho \right) $, the principal

%Since at each point $\varphi _k$ the derivative with respect to $\varphi $
%of the left-hand side of equation (\ref{eq8}) does not vanish for $\rho =0$,

% such that $\varphi _k\left( 0\right) =\varphi _k$.
%This unique solution determines the admissible (with respect  to  a
%small loop  encircling the point $z_0$) cut
%near the point $z_0$. So, the boundary condition for an admissible cut near
%the singular points is%
%$$
%\lim \limits_{\rho \rightarrow 0}\varphi \left( \rho \right) =\varphi _k,
%$$
%where $\varphi _k$ are given by formula (\ref{eq9}). 

%The above consideration can be summarized in the form of the following theorem.

%\newline

Remark that 
 if  $\dot z_a = 0$  
%that is, $z_0$ is a stationary singular point, 
and $\dot\Phi\left( z_a(t),t\right)\ne 0 $, from formula \eqref{strW} follows
$$
{\rm var}_l\,p_j =-\frac{2}{3k_j}\rho ^{\frac {3}{2}}\, \Re \,%
\,\left[ e^{\left( \frac{3i\varphi }{2}+i\theta _0\right) } \{\dot R_0+iR_0\dot\theta _0\}  \left( 1+\xi _1 \left(
z_a+\rho e^{i\varphi }\right) \right)     \right] =0,
$$
%and  the directions of admissible cuts near the singular point are defined by the equation
%\begin{equation}\label{eq81}
%\dot R_0 \cos \left( \frac{3\varphi }{2}+\theta _0\right)- R_0\dot\theta _0\sin \left( \frac{%
%3\varphi }{2}+\theta _0\right)  =0,
%\end{equation}
%which can be rewritten in the form
%$
%\sin \left( \frac{3\varphi }2+\theta _0-\nu _0\right) =0,
%$
which results in three direction of admissible cuts
$
\varphi =\varphi _k = \frac{2}{3}(\pi k-\theta _0 +\nu _0)$,
$k=0, \pm 1, \pm 2,...$,
where $\nu _0=\arcsin (\dot R_0/\sqrt{\dot R_0^2+R_0^2\dot \theta _0^2)}$, 
$R_0$ and $\theta _0$ are  the modulus and the argument of $\Phi\left( z_a(t),t\right) $ respectively.
In the case when  $\dot z_a = 0$ and first $(j-1)$ time derivatives of $\Phi\left( z,t\right) $ at $z_a(t)$ equal to zero, but the $j$-th derivative is not, the number of directions is $(2j+3)$.

%In the case $(c)$, when $\dot z_0=\dot c_0=...=\dot c_{j-1}=0$ and $\dot c_j\ne 0$,
%$
%\varphi =\varphi _k = \frac{2\pi k}{2j+3}-\theta _j +\nu _j,$
%$k=0, \pm 1, \pm 2,... ,
%$
%where $\nu _j=\arcsin (\dot R_j/\sqrt{\dot R_j^2+R_j^2\dot \theta _j^2)}$ and $c_j=R_je^{i\theta _j}$.

The constructed support of the distribution $\mu _j$ must satisfy the following 
 conditions: each cut emanates from a singular point of $W_j$,  located in the domain $\Omega _j (t)$, and  ${\rm var}_l\,p_j$ vanishes on each cut. To obtain a  two-phase mother body, one has to calculate the  corresponding density along each cut. 
%Note that, here we do not discuss the existence of the two-phase mother body in general situation.
In the next section we show examples when a two-phase mother body exists and is unique, and we use the constructed mother bodies to derive the exact solutions to the Muskat problems.
 
%\section{ Algorithm }\label{sec:alg}

%Assume that the droplet $\Omega(0)$ has an algebraic boundary and the singularities of the complex potential $W$, including the singularity at infinity, are not more than those of logarithmic type.

%The suggested algorithm consists of the following steps.
%\newline
%{\it Step 1} Determine the singularities of the Schwarz function of $\Omega (0)$ by solving an algebraic equation defining the boundary in terms of characteristic variables $(z,w)$. 
%\newline
%{\it Step 2} Construct a set of admissible cuts.
%\newline
%{\it Step 3} Select the set of bounded admissible cuts.
%\newline
%{\it Step 4} Calculate a dynamical distribution of the sinks along the cuts.

\section{ Examples of specific $\Gamma (0)$ \label{sec:ex}}

\subsection{Circle } 

To illustrate the method, we start with the simplest example for which solution is known. Suppose that the initial shape of the interface is a circle with equation $x^2+y^2=a^2(0)$, and during the evolution the boundary remains circular, 
 $x^2+y^2=a^2(t)$. The corresponding Schwarz function is $S=a^2(t)/z$. 
Equation \eqref{main} in this case reads as
$$\partial _z(W_2-W_1)=\Bigl(\frac{1}{k_1}-\frac{1}{k_2}\Bigr )a\dot a/z.$$
Integrating this equation, we have
$$W_2-W_1=\Bigl(\frac{1}{k_1}-\frac{1}{k_2}\Bigr )a\dot a \log z +C(t).$$
The complex potential $W_2$ has a singularity at zero, while $W_1$ has a singularity at infinity. Thus, the two-phase mother body has support at these two points, one of which serves as a sink and the other as a source. 
Taking the real parts of both sides of the previous equation, we obtain
$$p_2-p_1=\Bigl(\frac{1}{k_1}-\frac{1}{k_2}\Bigr )a\dot a \ln \sqrt{x^2+y^2} +\Re C(t),$$
which is satisfied if
$p_j=-\frac{a\dot a}{2k_j} \ln (x^2+y^2) + C_j(t)$, $j=1,2$ with $C_j$ chosen from the condition \eqref{2}. This choice is, obviously, not unique. To specify, for instance, 
$C_1(t)$ one could use the condition at infinity. If the condition reads as
$$p_1(x,y,t)=-\frac{Q(t)}{2\pi k_1}\ln \sqrt{x^2+y^2}\quad\mbox{as}\quad \sqrt{x^2+y^2}\to\infty$$
with a defined sink/source strength $Q(t)$ at infinity, then $C_1=0$ and  
$p_1=-\frac{a\dot a}{2k_1} \ln (x^2+y^2)$. The evolution of the boundary is defined by the equation $Q(t)=\dot A=2\pi a\dot a$ via $\dot a$. Here $\dot A(t)$ is the  rate of change of the area of the interior domain $\Omega _2 (t)$.
The interior pressure is
$p_2=-\frac{a\dot a}{2k_2} \ln (x^2+y^2) + \frac{a\dot a}{2k_2}\ln a^2-\frac{a\dot a}{2k_1}\ln a^2$ with a  source/sink  at the origin of the strength $|Q(t)|$. 

Alternatively, we could choose  $C_j$ such that $p_j$ vanishes on $\Gamma$,
\begin{equation}\label{circlep}
p_1=-\frac{a\dot a}{2k_1} \ln (x^2+y^2)+\frac{a\dot a}{2k_1}\ln a^2, \quad
p_2=-\frac{a\dot a}{2k_2} \ln (x^2+y^2) + \frac{a\dot a}{2k_2}\ln a^2.
\end{equation}
We remark, that in the case of the circular initial interface,  condition \eqref{2} could be replaced with $p_1-p_2=\gamma\kappa$  while keeping the boundary in the family of concentric circles during the course of the interface evolution. Here $\gamma$ is a constant surface tension coefficient and $\kappa$ is a free boundary curvature. In that case the first equation in \eqref{circlep} is replaced with
$p_1=\frac{a\dot a}{2k_1} \ln (a^2/(x^2+y^2))+\gamma /a$.

\subsection{Ellipse}

Consider a two-phase problem with an elliptical interface,
$\Gamma (0)=\left\{ \frac{x^2}{a(0)^2}+\frac{y^2}{b(0)^2}=1\right\}$,
  where $a(0)$ and $b(0)$ are given and $a(0)>b(0)$.
% (see Fig. \ref{fgElli}).
The Schwarz function of an elliptical interface with  semi-axes $a(t)$ and $b(t)$ is
$$
S\left( z,t\right)
=\Bigl ( \bigl ( a(t)^2+b(t)^2 \bigr ) z-2a(t)b(t)\sqrt{z^2-d(t)^2}\Bigr )/d(t)^2,
$$
where $d(t)=\sqrt{a(t)^2-b(t)^2}\,$ is the half of the inter-focal distance.
Assuming that 
the interface remains elliptical during the course of the evolution, 
%the latter is justified by the computations below, 
from equation (\ref{main}) we have
%\begin{equation}\label{wzellipse}
%\pd{W}{z}=-\frac{z}{2}\pd{}{t}\Bigl (  \frac{a^2+b^2}{d^2}   \Bigr )
%+\sqrt{z^2-d^2}\,\pd{}{t}\Bigl (  \frac{ab}{d^2}   \Bigr )
%-\frac{1}{ \sqrt{z^2-d^2}}\, \frac{ab}{2d^2}\pd{}{t}\Bigl (  d^2   \Bigr )
%\end{equation}
%and
\begin{eqnarray}
W_2-W_1=\Bigl(\frac{1}{k_2}-\frac{1}{k_1}\Bigr )\Bigl \{-\frac{z^2}{4}\pd{}{t}\Bigl (  \frac{a^2+b^2}{d^2}   \Bigr ) \notag \\
+\frac{z}{2}\sqrt{z^2-d^2}\,\pd{}{t}\Bigl (  \frac{ab}{d^2}   \Bigr )
- \frac{1}{2}\log (z+\sqrt{z^2-d^2})\, \pd{(ab)}{t}\Bigr \}+C(t).\label{well}
\end{eqnarray}
%Here $p=\Re [W]$, and the stream function $\psi =-\Im [W]$.
%\begin{figure}
%\hskip1.5cm
%\begin{center}
%\includegraphics[width=0.3\textwidth]{PicturesHeSh1_Elli1.eps}{a}
%\includegraphics[width=0.3\textwidth]{PicturesHeSh1_Elli2.eps}{b}
%\includegraphics[width=0.3\textwidth]{PicturesHeSh1_Elli3.eps}{c}
%\caption{Ellipse: (a) $t=0$, $a=2$, $b=1$, $d^2=3$; (b) first scenario for $t>0$ with %$d=const$;
%(c) second scenario $ab/d^2=const$. } 
%\label{fgElli}
%\end{center}
%\end{figure}

The first two terms in the right hand side of (\ref{well}) have  poles of order two at infinity.
%The complex potential $W$ cannot have singularities in the finite domain  $\mathbb{R}^2\setminus\Omega (t)$.
Those terms are  eliminated, if the eccentricity of the ellipse does not change with time, the latter implies that the ratio $a(t)/b(t)=const$.
% Thus,    can grow for all times only if the second term in the right hand side of \eqref{well} is eliminated, which
This ensures the existence of no more than logarithmic singularity at infinity and agrees with the solution to the exterior one-phase problem reported in \cite{HowBubble}. Thus, the expression for the complex potentials reduces to
\begin{equation}
W_2-W_1=
 -\frac{1}{2}\pd{(ab)}{t}\Bigl(\frac{1}{k_2}-\frac{1}{k_1}\Bigr )\log (z+\sqrt{z^2-d^2}) +C(t).\label{wellm}
\end{equation}
%$W_1-W_2=-\frac{1}{2}\pd{ab}{t}\Bigl(\frac{1}{k_1}-\frac{1}{k_2}\Bigr )\log (z+\sqrt{z^2-d^2}) + C(t)$.
Taking the real parts of both sides, we have 
$$p_2-p_1=
 -\frac{1}{2}\pd{(ab)}{t}\Bigl(\frac{1}{k_2}-\frac{1}{k_1}\Bigr )\ln |z+\sqrt{z^2-d^2}| +C_2(t)-C_1(t),
$$
where $C_j(t)$ may be chosen from the condition $p_j=0$ on $\Gamma (t)$, which leads to 
\begin{equation}
p_j=
 -\frac{1}{2k_j}\pd{(ab)}{t}\Bigl (\ln |z+\sqrt{z^2-d^2}| -\ln (a+b)\Bigr ),\label{wellm1}
\end{equation}
or
$$
p_j=
 -\frac{1}{2k_j}\pd{(ab)}{t}\Bigl (\ln \sqrt{(x+\alpha )^2 (1+y^2/\alpha ^2) } -\ln (a+b)\Bigr ),
$$
where 
$$
\alpha ^2=\Bigl ( x^2-y^2-d^2 +\sqrt{(x^2-y^2-d^2)^2+4x^2y^2}  \Bigr ) /2.
$$
%where $\Re [C_1(t)]=-\frac{1}{2}\pd{(ab)}{t}\ln (a+b)$.
Note that the inter-focal distance, $d(t)=2b(t)\sqrt{a^2(0)/b^2(0)-1}$, of such an ellipse changes, while the eccentricity is constant. 
The  support of the two-phase mother body  consists of a point sink/source  at infinity and a  source/sink
 distribution with density $\mu _2(x,t)=2ab/(d^2 k_2) \, \partial _t(\sqrt{d^2-x^2})$
along the inter-focal segment.
%, $Q_2(t)/k_2=\pm \int _{-d}^{d} q_2(x,t)\, dx= \pm\pi\partial _t(ab)/k_2$.
%as well as (b) when the inter-focal distance is constant and the distribution $q(x,t)=2\partial _t(ab)\sqrt{d^2-x^2}$ \cite{dmb}.
%We remark that  to find the directions of the cut at the singular points of the Schwarz function, we used the formula \eqref{dir}, which 
The support of this distribution is defined using formula \eqref{dir}. 
Indeed, the singular points of $W_2$ are   $z=\pm d$
with $\Phi (z,t)=-2ab/d^2\, \sqrt{z\pm d}$ respectively.
 Formula \eqref{dir} implies that the direction of the cut 
at $z=d$ is defined by the angle
$\varphi =\pi +2\pi k$,  and at $z=-d$ by the angle $\varphi =2\pi k$, $k=0,\pm 1,\pm 2,\dots$.
Thus,  the two-phase mother body, described above, allows the interface between two fluids remain elliptical 
for an infinite time if the domain $\Omega _2$  grows (sources are located along the interfocal segment and a sink is located at the point of infinity). The opposite sink/source choice  allows the complete removal of the fluid initially occupied domain $\Omega _2$.

We remark that from (\ref{wellm1}) follows that 
 the pressure at infinity grows as $p_1 \sim -(2k_1)^{-1}\ln|z| \partial _t(ab)=-\frac{\dot A}{2k_1}\ln \sqrt{x^2+y^2}$, which agrees with formula \eqref{4m}.
Moreover, the strength of the sink/source at infinity is in agreement with the total strength of the source/sink distribution in $\Omega _2$ since $\int _{-d}^{d} k_2 \mu _2(x,t)\, dx= \pi\partial _t(ab)=\dot A$.

%implies a relation between the strength of the sink $Q_1(t)$ at infinity and the change of the area,   
%$$
%Q_1(t)=\pi \frac{d \bigl (a(t)b(t)\bigr)}{dt}=\frac{2\pi a(0)b(t)\dot b(t)}{b(0)}=\frac{2\pi b(0)a(t)\dot a(t)}{a(0)}.
%$$ 
%Indeed, $Q(t)=\partial _t A(t)/k_1$, where $A(t)$ is the area of $\Omega _2 (t)$.
%For ellipse $A(t)=\pi a(t)b(t)$, therefore, $Q_1(t)=\pi \partial _t (ab)$.   

%For the specific case when $Q=dA(t)/dt\equiv const$,  the bubble area grows linearly in time.

Note that Crowdy \cite{crowdy2006}
 obtained an exact solution with a different type of growth at infinity; the solution, reported in \cite{crowdy2006}, has a linear far field flow and
 a constant area of ellipse.
 We observe that under the assumption that the area of the elliptical inclusion  does not change in time, that is, $a(t)b(t)=const$, equation \eqref{well} implies
$$
W_j=\frac{1}{k_j}\Bigl \{-\frac{z^2}{4}\pd{}{t}\Bigl (  \frac{a^2+b^2}{d^2}   \Bigr ) 
+\frac{z}{2}\sqrt{z^2-d^2}\,\pd{}{t}\Bigl (  \frac{ab}{d^2}   \Bigr )
\Bigr \}+C_j(t), \quad j=1,2,
$$
therefore,  the pressure $p_1$ is defined by
$$
p_1=\frac{1}{k_1}\Re \Bigl \{-\frac{z^2}{4}\pd{}{t}\Bigl (  \frac{a^2+b^2}{d^2}   \Bigr ) 
+\frac{ab\,z}{2}\sqrt{z^2-d^2}\,\pd{}{t}\Bigl (  \frac{1}{d^2}   \Bigr )
\Bigr \} - \frac{b^2a \dot a}{k_1 d^2},
$$
which retries a linear far field flow. However, the obtained solution is different from Crowdy's, since the interior flow reported in \cite{crowdy2006} is a simple linear flow, while the solution in question, in addition to the linear flow, has another term
$$
p_2=\frac{1}{2k_2}\Bigl \{\frac{(y^2-x^2)}{2}\pd{}{t}\Bigl (  \frac{a^2+b^2}{d^2}   \Bigr ) 
+\frac{ab\,x(\alpha ^2-y^2)}{\alpha}\,\pd{}{t}\Bigl (  \frac{1}{d^2}   \Bigr )
\Bigr \} - \frac{b^2a \dot a}{k_2 d^2}.
$$
The interior flow is generated by the density 
$$
\mu = \frac{ab\,\partial _t (d^2)}{k_2d^4}\,\,\frac{(2x^2-d^2)}{\sqrt{d^2-x^2}},
$$
supported on the inter-focal segment. Such a density changes sign along the inter-focal segment, so the area of the ellipse does not change in time: if $a(t)$ increases with time, the ellipse becomes ``thiner''.

\subsection{The Neumann's oval}

Let the initial free boundary have a shape of the Neumann's oval  \cite{shapiro} given by the equation
$\Gamma (0) =\left\{   (x^2+y^2)^2 -a(0)^2x^2-b(0)^2y^2 =0\right\}$, (see Fig. \ref{fgNeu}). Its Schwarz function is
$$
S(z,0)=\Bigl (z(a(0)^2+b(0)^2)+ 2z\sqrt{z^2d(0)^2+a(0)^2b(0)^2}\Bigr )/(4z^2-d(0)^2),
$$
where $d(0)^2=a(0)^2-b(0)^2>0$ with given $a(0)$, $b(0)$.
Assume that during the evolution the domain retains the Neumann's oval shape,
$$\Gamma (t) =\left\{   (x^2+y^2)^2 -a(t)^2x^2-b(t)^2y^2 =0\right\},$$
 with unknown $a(t)$, $b(t)$ for $t>0$. The singularities of the Schwarz function, 
$
S(z,t)=\Bigl (z(a(t)^2+b(t)^2)+ 2z\sqrt{z^2d(t)^2+a(t)^2b(t)^2}\Bigr )/(4z^2-d(t)^2),
$
 located in the interior domain $\Omega _2 (t)$ are simple poles at $z=\pm d/2$, while the singularities located
in the exterior domain $\Omega _1$ are  the branch points at  $z=\pm ia(t)b(t)/d(t)$.

To ensure at most logarithmic growth of the pressure at the singular points,  $d$ must be constant. In that case
\begin{equation}
\dot S =\frac{z\partial _t (a^2+b^2)}{4z^2-d^2}+ \frac{z\partial _t (a^2b^2)}{(4z^2-d^2)\sqrt{z^2d^2+a^2b^2}}.
\end{equation}
Then equation \eqref{main} implies
\begin{eqnarray}\label{mainNeu}
{p_2}- {p_1}=-\frac{1}{2}\Bigl(\frac{1}{k_2}-\frac{1}{k_1}\Bigr )\notag
\Bigl \{\frac{\partial _t(a^2+b^2)}{8}\Re \bigl [\log (4z^2-d^2)\bigr ] \\-
\frac{\partial _t(a^2b^2)}{2(a^2+b^2)} \Re \bigl [\tanh ^{-1} \frac{2\sqrt{a^2b^2+d^2z^2}}{a^2+b^2}     \bigr ]    \Bigr \}+C(t).
\end{eqnarray}
Taking into account that when $d$ is constant, $\partial _t (a^2+b^2)=4a\dot a$ and $\partial _t (a^2b^2)=2a\dot a (a^2+b^2)$, we have
\begin{eqnarray}\label{mainNeu1}
{p_2}- {p_1}=-\frac{a\dot a}{2}\Bigl(\frac{1}{k_2}-\frac{1}{k_1}\Bigr )\notag
\Re\Bigl \{  \log (4z^2-d^2) \\-
\log (a^2+b^2+2\sqrt{a^2b^2+d^2z^2}\,)\Bigr \}+\ln(d)         +C(t).
\end{eqnarray}
Note that on the interface 
$$
\log \frac{4z^2-d^2}{a^2+b^2+2\sqrt{a^2b^2+d^2z^2}}=\log \frac{z}{\bar z},
$$
whose real part is zero on $\Gamma (t)$. Therefore, to satisfy the condition \eqref{2},
$$
C(t)=-\frac{a\dot a \ln d}{2}\Bigl(\frac{1}{k_2}-\frac{1}{k_1}\Bigr ).
$$
Thus, 
we have 
\begin{equation}\label{pneu}
p_j=\frac{a\dot a}{2k_j}\Re \log \frac{a^2+b^2+2\sqrt{a^2b^2+d^2z^2}}{4z^2-d^2},\qquad j=1,2,
\end{equation}
and the interior part of mother body $\mu _2$ consists of either two point sinks or two point sources.

To find the  directions of the cuts in $\Omega _1$, we use formula \eqref{dir} (it is  general position like an ellipse). In the neighborhood of the branch point $z_0=iab/d\in \Omega _1$, $\arg [\Phi\left( z_0(t),t\right)]=-\pi/4+\pi k$,     $\arg [\dot z_0]=\pm \pi/2$, where the plus corresponds to the growth of the interior domain $\Omega _2$. The direction of the cut near this point is $\varphi =\pi /2+ 2\pi k$, $k=0,\pm 1,\pm 2, \dots$.
Similarly, at the $z_0=-iab/d$, $\arg [\Phi\left( z_0(t),t\right)]=\pi/4+\pi k$,     $\arg [\dot z_0]=\pm \pi/2$, where the plus corresponds to the decrease of the interior domain $\Omega _2$. The direction of the cut near this point is $\varphi =-\pi /2+ 2\pi k$, $k=0,\pm 1,\pm 2, \dots$.
%$$W=\frac{i}{2}d_t[abd/(a^2+b^2)]Z \mp
%  (1/3)d_t[(2ab)^{3/2}/d^{1/2}](iZ)^{3/2} \mp 
%(iZ/d)^{1/2}(2ab)^{3/2}d_t(ab/d)+O(Z). $$
\begin{figure}
\hskip1.5cm
\begin{center}
\includegraphics[width=0.3\textwidth]{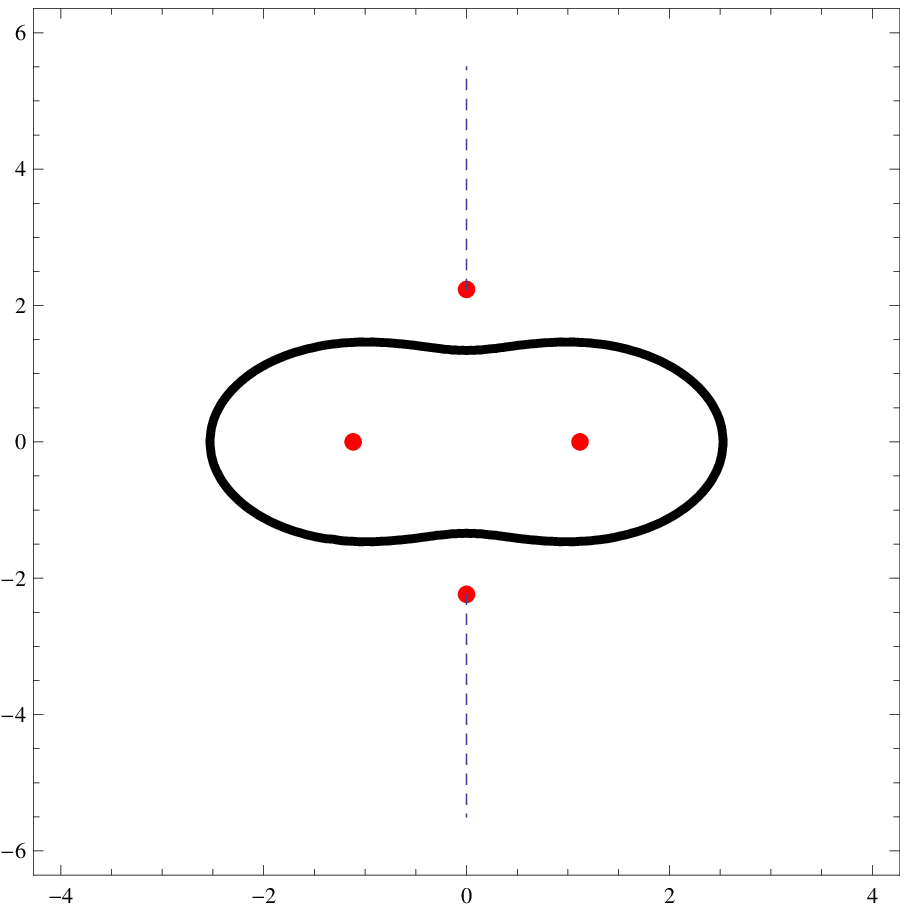}{a}
\includegraphics[width=0.3\textwidth]{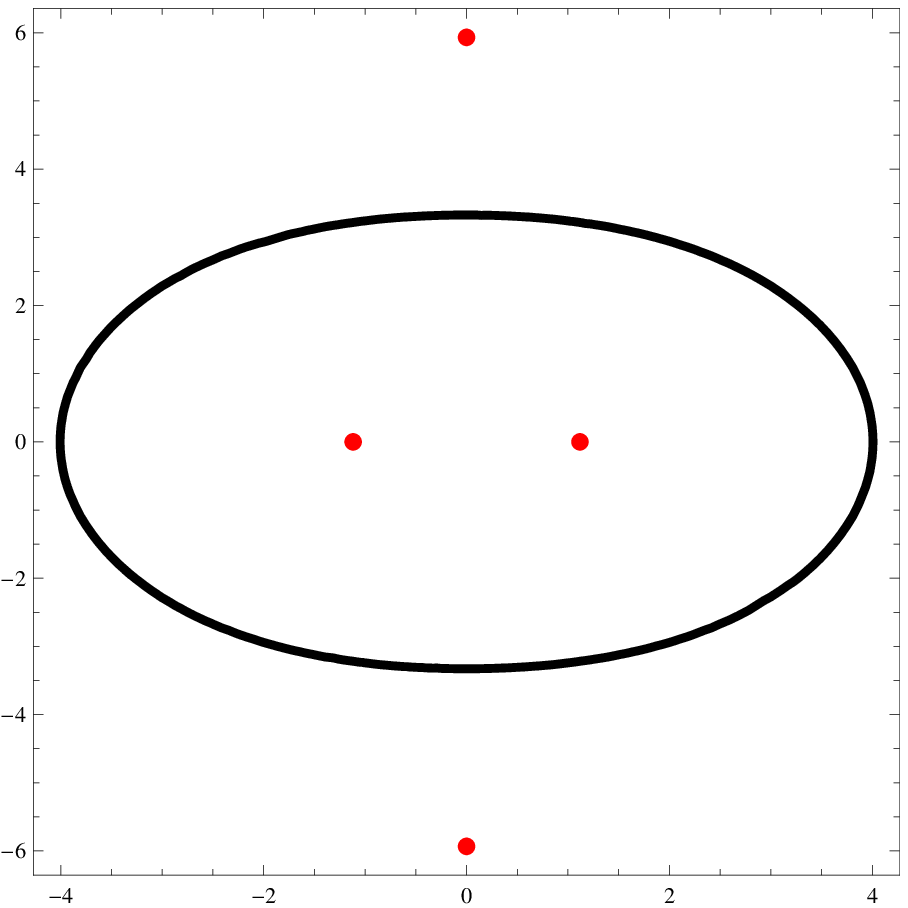}{b}
\caption{The Neumann's ovals (solid lines), singularities of the complex potential (dots), and the cuts (dashed lines) for  $d=\sqrt{5}$: (a) $a=2.5$, $b=\sqrt{5}/2$;  (b) $a=4$, $b=\sqrt{11}$. 
%(c) $a=3$, $b=2$.
} \label{fgNeu}
\label{fgNeu}
\end{center}
\end{figure}
The support of the two-phase mother body for the Neumann's oval is shown in Fig.  \ref{fgNeu}. The cuts are the dashed lines that go along the imaginary axis starting at each branch point (the dots in the exterior domain in Fig.  \ref{fgNeu}) to infinity. The dots in the interior domain correspond to the simple poles.

%To ensure at most logarithmic growth  at infinity, one should choose $d$ to be constant. Indeed, then the leading term of $\partial _z W$ at infinity is
%$$\partial _zW \sim -\frac{\partial _t(a^2+b^2)}{8z}. $$
%The latter expression defines  
%Note that $v_y$ is continuous along the cuts, and the desired sink distribution $q(y,t)$ can be computed via a jump of $v_x$. 
To obtain the sink/source density along the cut located above the $x$-axis, we first compute the variation of $S(z,t)$
$${\rm var}_l \, S(z)={4z(4z^2-d^2)^{-1}\sqrt{z^2d^2+a^2b^2}}{\bigl |_{z=iy, y>ab/d}}=
\frac{4y\sqrt{y^2d^2-a^2b^2}}{4y^2+d^2} ,$$
then the jump of $\partial _zW_1$, that  is,
$${\rm var}_l\, \partial_zW_1=-\frac{1}{2k_1}\partial_t({\rm var}_l \, S(z))=
-\frac{2y}{k_1}\partial_t \Bigl (\frac{\sqrt{y^2d^2-a^2b^2}}{4y^2+d^2}\Bigr ), $$
finally,  the sink distribution on both cuts equals
 $$
\mu_1(y,t)=\frac{1}{k_1}\Bigl | \frac{\partial _t(a^2b^2) y}{(4y^2+d^2)\sqrt{y^2d^2-a^2b^2}}\Bigr |. 
$$

To compute the  rate, $Q_1(t)$, through the cuts located in the exterior domain, one has to integrate  $k_1\mu_1(y,t)$ along the cuts, which implies
\begin{equation}\label{420}
Q_1(t)=\frac{\pi}{2}\frac{\partial _t (a^2b^2)}{(a^2+b^2)}=\frac{\pi}{2}\partial _t (a^2).
\end{equation}

From formula \eqref{pneu} follows that the  rate at infinity is
\begin{equation}\label{419}
Q(t)=\frac{\pi}{4}\partial _t (a^2+b^2)=\frac{\pi}{2}\partial _t (a^2).
\end{equation}

Those rates are linked to the change of the area
\begin{equation}\label{421}
 Q(t)+Q_1(t) =\pi \partial _t (a^2)=\dot A,
\end{equation}
where $A(t)$ is the area of the interior domain.
% and dot denotes the time derivative.
Indeed, the area of the Neumann's oval is $A=\pi(a^2+b^2)/2$ (see \cite{shapiro}, p. 20), which can be rewritten as $A=\pi (a^2-d^2/2)$.
Since in the case in question $d$ is constant,  $\dot A =\pi\partial _t(a^2)$.

Note also that
the rate through each of the two point sources/sinks located in the domain $\Omega _2$  equals
$\dot A/2$.

We remark that since $d^2(t)=a^2(t)-b^2(t)=const$, in the case of decreasing area of $\Omega _2$, the obtained solution is valid up to  the limiting case when $b$ approaches zero and $\Gamma$ splits into two circles $(x+\frac{d}{2})^2+y^2=\frac{d^2}{4}$ and 
$(x-\frac{d}{2})^2+y^2=\frac{d^2}{4}$.

\subsection{The  Cassini's oval}

\begin{figure}
\hskip1.5cm
\begin{center}
\includegraphics[width=0.3\textwidth]{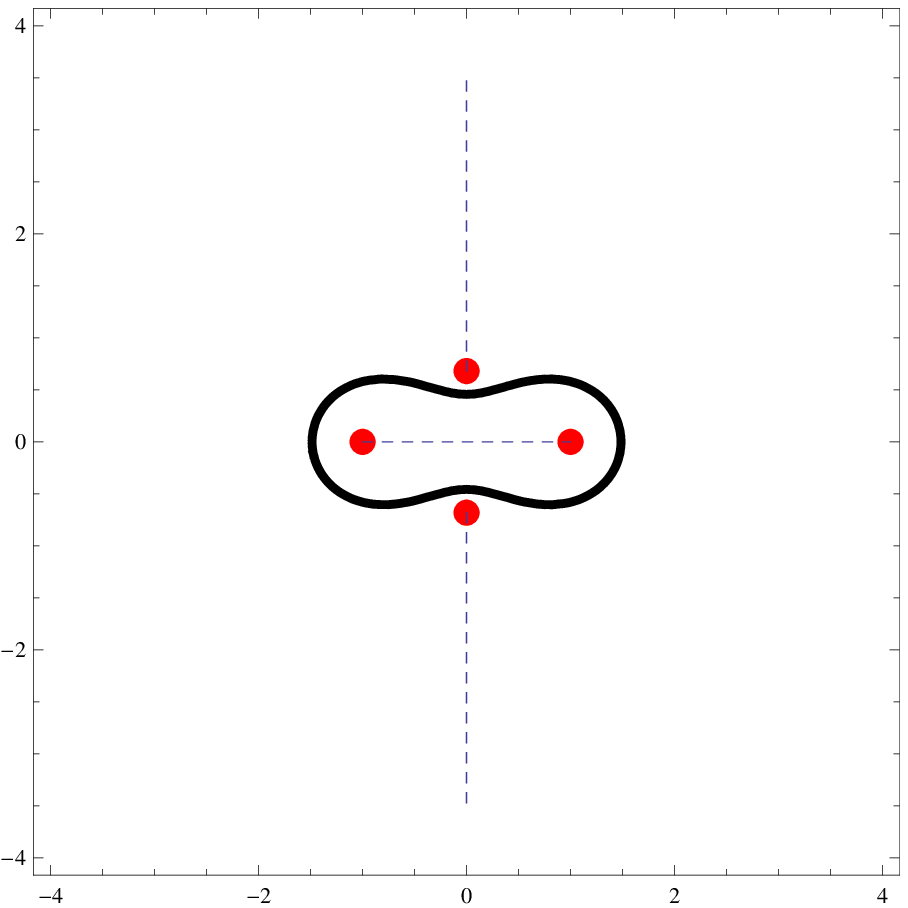}{a}
\includegraphics[width=0.3\textwidth]{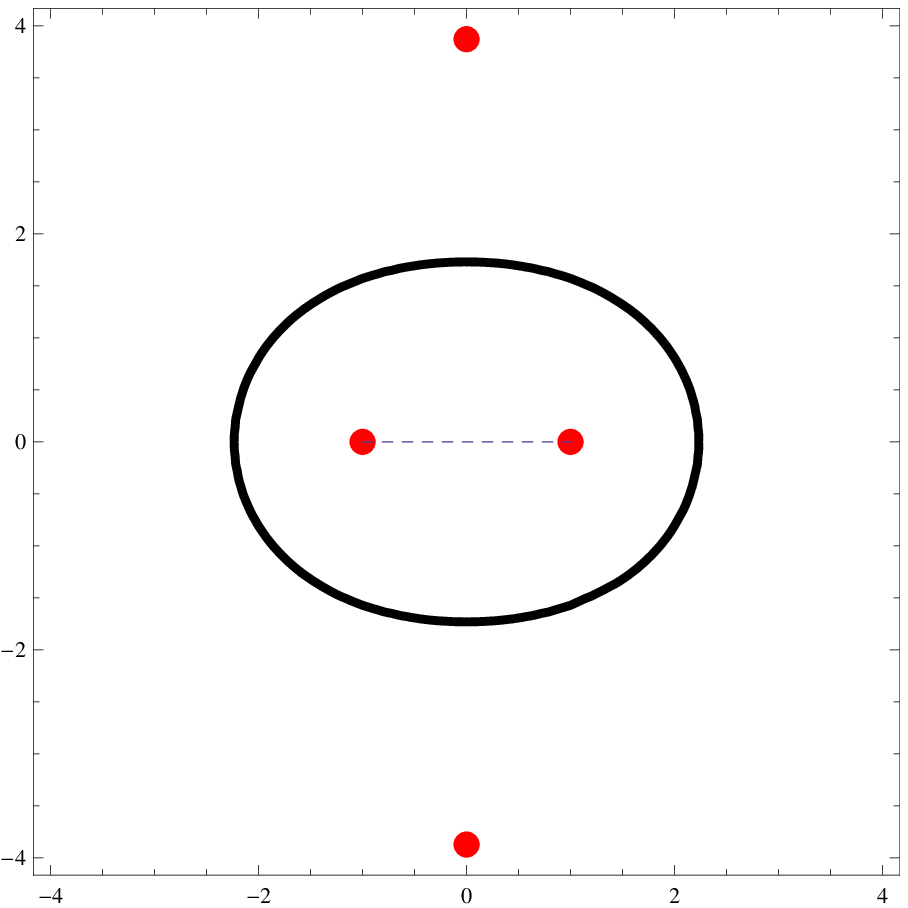}{b}
\caption{The Cassini's ovals  (solid lines), singular points of the complex potential (dots), and the cuts (dashed lines) for $b=1$: (a) $a=1.1$, 
%(c) $a(t_1)=1.5$, 
 (b) $a=2$. 
%The dots correspond to the singularities, and the dashed lines correspond to the cuts.
} 
\label{fgCass}
\end{center}
\end{figure}

Similar to the previous examples, assume that $\Gamma (t)$   remains in the specific family of curves,
the Cassini's ovals, given by the equation
$$
\left( x^2+y^2\right) ^2-2b(t)^2\left( x^2-y^2\right) =a(t)^4-b(t)^4,
$$
%or
%$$
%\left [( x-b)^2+y^2\right] \left[( x+b)^2+y^2\right ] =a^4,
%$$
where $a(t)$ and $b(t)$ are unknown positive functions of time.
% (see Fig. \ref{fgCass}).
This curve consists  of one
closed curve if $a(t)>b(t)$ (see Fig. \ref{fgCass}), and  two closed curves otherwise. Assume that at $t=0$ $a(0)>b(0)$. 
%This example is interesting from the following point of view. If we
%compute 
The Schwarz function of the Cassini's oval,
$$
S\left( z,t \right) =\sqrt{b^2z^2+a^4-b^4}\, /\sqrt{z^2-b^2},
$$
% has  two singular points $z=\pm b$ inside the domain $\Omega $ with singularities of  the  %type  of  the inverse square root, which is not the generic case.
%The domain and the singular  points  are shown in Figure~\ref{fg10}. 
 has two singularities in $\Omega _1(t)$, $z=\pm i\sqrt{(a^4-b^4)/b^2}$, and two singularities, $z=\pm b$, in $\Omega _2(t)$.
%This happens since these singularities are
%generated not by regular characteristic points of the surface
%$\Gamma _{\mathbb{C}}$ but by singular
%points of this surface.%
%It is known that for $a=b$, when the oval looks like eight-curve, area of each loop is $a^2$:
%\begin{equation}
%\int\limits _0 ^a \frac{ax\,dx}{\sqrt{a^2-x^2}}=-a\sqrt{a^2-x^2}|_0 ^a =a^2.
%\end{equation} 
To ensure that the singularities of the complex potential have no more than the logarithmic type, $\dot b$ must be zero.  Thus, we have
$$
\dot S\left( z\right) =\frac{2a^3\dot a}{\sqrt{b^2z^2+a^4-b^4}\sqrt{z^2-b^2}},
$$
%and
%$
%\pd{W}{z} =-\frac{a^3\dot a}{\sqrt{b^2z^2+a^4-b^4}\sqrt{z^2-b^2}}
%$
%As long as $a(t)>1$, the inner singularities are at points $x=\pm 1$ only and are integrable.
Then equation \eqref{main} implies
\begin{equation}\label{mainCas}
{W_2}- {W_1}=-a\dot a\Bigl(\frac{1}{k_2}-\frac{1}{k_1}\Bigr )
\,F \Bigl (\cos ^{-1}  \bigl (\frac{b}{z} \bigr),\, \frac{\sqrt{a^4-b^4}}{a^2}\,\Bigr )+C(t),
\end{equation}
%and $W=-a\dot a\,F(\arccos \, b/z,\, \sqrt{a^4-b^4}\, /a^2)+C(t)$,
 where $F(\alpha\,,\beta)$ is the incomplete elliptic integral of the first kind,
$$
F \Bigl (\cos ^{-1} \bigl ( \frac{b}{z}\bigr ),\, \frac{\sqrt{a^4-b^4}}{a^2}\, \Bigr )=\int\limits _0 ^{\sqrt{1-b^2/z^2}}
\frac{dt}{\sqrt{1-\frac{a^4-b^4}{a^4}t^2}\,\sqrt{1-t^2}}=
\int\limits _0 ^{\cos ^{-1}  (b/z )}\frac{dt}{\sqrt{1-\frac{a^4-b^4}{a^4}\sin ^2 t}}.
$$
Then
\begin{equation}\label{mainCas1}
{p_2}- {p_1}=-\frac{a\dot a}{2}\Bigl(\frac{1}{k_2}-\frac{1}{k_1}\Bigr )
\,\Bigl [F \Bigl (\xi,\, \frac{\sqrt{a^4-b^4}}{a^2}\,\Bigr )+\overline{F \Bigl (\xi,\, \frac{\sqrt{a^4-b^4}}{a^2}}\,\Bigr )\Bigr]+C(t),
\end{equation}
where $\xi=\cos ^{-1}  \bigl (\frac{b}{z} \bigr)$.
Using the property $\overline {F(\alpha\,,\beta)}=F(\overline{\alpha}\,,\beta)$ and the summation formula for the elliptic integrals \cite{BE}, we have
\begin{equation}\label{mainCas1m}
{p_2}- {p_1}=-\frac{a\dot a}{2}\Bigl(\frac{1}{k_2}-\frac{1}{k_1}\Bigr )
\,F \Bigl (\alpha
,\, \frac{\sqrt{a^4-b^4}}{a^2}\,\Bigr )+C(t),
\end{equation}
where
\begin{equation}\label{alpha1}
\alpha=\sin ^{-1}\frac{
\cos \overline{\xi}\sin \xi \sqrt{1-\frac{{a^4-b^4}}{a^4}\sin ^2\overline{\xi}}+
\cos \xi\sin \overline{\xi} \sqrt{1-\frac{{a^4-b^4}}{a^4}\sin ^2{\xi}}
}
{   1-\frac{{a^4-b^4}}{a^4}\sin ^2 \xi \sin ^2 \overline{\xi}   }.
\end{equation}
Rewriting \eqref{alpha1} in terms of $z$ and $\overline z$, we obtain the following expression
%%%%%%%%%%%%%%%%%%%%%%%%%% Intermediate expression for alpha if needed
%\begin{equation}\label{alpha2}
%\alpha=\sin ^{-1}\frac{
%\frac{b}{\bar z}\sqrt{1-\frac{b^2}{z^2}} \sqrt{1-\frac{{a^4-b^4}}{a^4}\frac{(\bar z^2-b^2)}{\bar z^2}}+
%\frac{b}{z}\sqrt{1-\frac{b^2}{\bar z^2}}\sqrt{1-\frac{{(a^4-b^4)}}{a^4}\frac{(z^2-b^2)}{z^2}}
%}
%{   1-\frac{{(a^4-b^4)}}{a^4}\frac{(z^2-b^2)}{z^2}  \frac{(\bar z^2-b^2)}{\bar z^2} },
%\end{equation}
%%%%%%%%%%%%%%%% Final expression for alpha
\begin{equation}\label{alpha3}
\alpha=\sin ^{-1}\frac{
a^2z\sqrt{z^2-b^2}\sqrt{b^2\bar z^2+a^4-b^4}+
a^2\bar z \sqrt{\bar z ^2 -b^2} \sqrt{b^2z^2+a^4-b^4}
}
{   b^2z^2\bar z ^2+(a^4-b^4)(z^2+\bar z ^2-b^2)  }.
\end{equation}
The pressures satisfying \eqref{mainCas1m} are
\begin{equation}\label{mainCas2}
{p_j}=-\frac{a\dot a}{2k_j}
%\Bigl(\frac{1}{k_2}-\frac{1}{k_1}\Bigr )
\,F \Bigl (\alpha
,\, \frac{\sqrt{a^4-b^4}}{a^2}\,\Bigr )+C_j(t),
\end{equation}
where
the terms $C_j(t)$ are computed from the values of $p_j$ on the interface, on which
$$
{p_j}=-\frac{a\dot a}{2k_j}\,F \Bigl (\frac{\pi}{2}
,\, \frac{\sqrt{a^4-b^4}}{a^2}\,\Bigr )+C_j(t). 
$$
Finally, for the pressure we have  
\begin{equation}\label{mainCas3}
{p_j}=-\frac{a\dot a}{2k_j}
%\Bigl(\frac{1}{k_2}-\frac{1}{k_1}\Bigr )
\,F \Bigl (\alpha
,\, \frac{\sqrt{a^4-b^4}}{a^2}\,\Bigr )+\frac{a\dot a}{2k_j}\,F \Bigl (\frac{\pi}{2}
,\, \frac{\sqrt{a^4-b^4}}{a^2}\,\Bigr ).
\end{equation}

Let us construct the two-phase mother body starting with its part located in the domain $\Omega _1$.
This is a generic situation, so we can use formula \eqref{dir}. 
In the neighborhood of the point $z_0= i\sqrt{(a^4-b^4)}\,/b$, $\arg [\dot z _0]=\pi /2 +2\pi k$,
$\arg  [\Phi\left( z_0(t),t\right)]=-\pi /4+\pi k$. Thus, according to  \eqref{dir}  the direction of the cut is $\varphi =\pi /2+2\pi k$, $k=0,\pm 1,\pm2, \dots$.

Similarly, at the point 

$z_0= -i\sqrt{(a^4-b^4)}\,/b$, $\arg [\dot z _0]=-\pi /2 +2\pi k$,
$\arg  [\Phi\left( z_0(t),t\right)]=-3\pi /4+\pi k$. 
\newpage
Therefore, the direction of the cut is $\varphi =-\pi /2+2\pi k$.

Taking into account symmetry with respect to $x$-axis, we conclude that the 
 support  of $\mu _1$ consists of  two rays starting at the branch points and going to infinity.  
The corresponding  density is
$$
\mu_1(y,t) =-\frac{2a^3\dot a}{k_1\sqrt{b^2y^2-a^4+b^4}\sqrt{y^2+b^2}}.
$$
%\begin{figure}
%\hskip1.5cm
%\begin{center}
%\includegraphics[width=0.4\textwidth]{PictHeShExtStream2.eps}
%%\includegraphics[width=0.3\textwidth]{PictHeShExtTre2.eps}{b}
%%\includegraphics[width=0.3\textwidth]{PictHeShExtTV3.eps}{c}
%\caption{Streamlines for the  Cassini's oval  with $a=10$ and $b=1$; the bold line in the center is the oval, the solid curves are the %streamlines, the dashed lines are cuts emanating from the singular points. } \label{fgStream}
%\end{center}
%\end{figure}

The singularities of the Schwarz function in the interior domain $\Omega _2$ have the inverse square root type (which is not the generic case),
\begin{equation}\label{insideDir1}
S\left( z,t\right) =\frac{\Phi\left( z,t\right)}{\sqrt{z-z_0}},
\end{equation}
where $\Phi\left( z,t\right)= \sqrt{b^2z^2+a^4-b^4}/\sqrt{z\mp b}$ is a regular functions of $z$ in the
neighborhood of the point $z_0=\pm b$ with $\Phi\left( z_0,t\right)\ne 0$.

Expanding   function
$\Phi\left( z,t\right) $  into  the Taylor series with respect to $z$ at the
point $(z_0,t)$, 
%say $\Phi(z,t)=\sum\limits _{j=0}^{\infty}c_j(t)\left( z-z_0(t)\right) ^j$. 
\begin{equation}
\label{eq677}
S\left( z,t\right) =\sum\limits _{m=0}^{\infty}c_m(t)\left( z-z_0\right) ^{m-1/2}.
\end{equation}
Differentiating   (\ref{eq677})  with respect to $t$, taking into account that $z_0$ is a stationary singularity since $\dot b=0$, we have:
\begin{equation}
 \dot S \left( z,t\right) =
\sum\limits_{m=0}^{\infty} \dot c_m \left( z-z_0\right) ^{m-1/2}
\end{equation}
Integration of the latter formula with respect to $z$ using (\ref{main2m}) implies:
$$
W_2 \left( z,t\right) =-\frac{\dot c_0}{2k_2}  \sqrt{z-z_0}\,\Bigl (1+\xi(z,t)\Bigr ),
$$
where $\xi  $ is a regular function near $z_0$ vanishing at
this point,
and $c_0=\Phi\left( z_0,t\right)$. 

The variation of the pressure along the loop $l$ is
\begin{equation}
{\rm var}_l\,p_2 =\Re \,\Bigl \{-\frac{\dot c_0}{k_2}  \sqrt{z-z_0}\,\Bigl (1+\xi(z,t)\Bigr ) \Bigr \}.
\end{equation}
Consider a small neighborhood of the point $z_0$, where  $z=z_0+\rho e^{i\varphi }$ with a small $\rho$.
Setting the principal part of ${\rm var}_l\,p_2$ to  zero, 
we have
\begin{equation}
{\rm var}_l\,p_2 =-\frac{|\dot c_0|}{k_2}\sqrt{\rho }\, \Re \,
\,\Bigl \{ \exp i\left ( \frac{\varphi }{2}+\pi n+\arg [\dot c_0] \right) \Bigr \}
%\left( 1+2\xi _1 \left( z_0+\rho e^{i\varphi }\right) \right) \Bigr \} 
=0,
\end{equation}
where $n=0,\pm 1,\pm 2, \dots$.

Thus, $\varphi =\pi -2\arg [\dot \Phi\left( z_0,t\right)]+2\pi k$, which implies that the support of $\mu _2$ is the segment $[-b,\,b]$ with the density
 $\mu_2=\frac{2a^3\dot a}{k_2\sqrt{b^2x^2+a^4-b^4}\sqrt{b^2-x^2}}$.
Integrating $k_j \mu_j$ along the corresponding cuts, one obtains the rate of change of the area of $\Omega _2$, which is given by the formula
%$A(t)=\int\limits _0^{\pi}\sqrt{a^4(t)-b^4(0)\sin ^2\theta}\, d\theta$  (\cite{willis}, p. 295), which implies the equation for $a(t)$
$$
\dot A=\partial _t(a^2)\,F(\pi,\frac{b^2}{a^2})=\\ \pi \, \partial _t(a^2)\, _2F_1 \,(\frac{1}{2},\frac{1}{2};1;\frac{b^4}{a^4}),
$$
where $_2F_1$
%$$
%F(m,n;c;z)=\sum\limits _{k=0}^{\infty}\frac{(m)_k(n)_k}{(c)_k}\frac{z^k}{k!}
%$$
is the hypergeometric series \cite{ww}.
% with $(m)_0=1$ and $(m)_k=m(m+1)\dots (m+k-1)$.
 %Fig. \ref{fgCass} shows three snapshots of an evolution of the Cassini's oval with $b(0)=1$ and $a(0)<a(t_1)<a(t_2)$, where $0<t_1<t_2$. 

%Note that for both cases, the Neumann's and Cassini's, when $a$ is large, and ovals have almost   circular shapes, their growth close to the boundary induces almost isotropic flow with a nonzero flux. In the case of the Neumann's oval, the sinks located on the external mother body slightly disturb this flow, while in the case of the Cassini's oval (see Fig. \ref{fgStream}), there is an intermediate region $|z|<<a^2$, where the flow is almost isotropic with $W_z \sim -\dot{a}/(az)$. However, at $|z| \sim a^2$ the flow changes, and for $|z|>>a^2$ the derivative of the complex potential $W_z \sim -\dot{a}a^3/z^2$, and sinks  suck all of the oil.

\section{Conclusions}\label{sec:concl}

We have studied  a Muskat problem with a negligible surface tension and suggested a method of finding exact solutions.   The idea of the  method was to keep the interface within a certain family of  curves defined by its initial shape by constructing two distributions with disjoint supports located on the different sides of the moving interface.

This study extended the results reported in \cite{howison2000} and \cite{crowdy2006}.
We gave new examples of exact solutions including the evolution of a circle, an ellipse, and two ovals: Neumann's and Cassini's. In those examples we assumed that the flux generated by the sinks/sources is finite, that is, the pressure may have at most a logarithmic growth. To demonstrate that this physical assumption does not restrict our method, we have presented an example of an exact solution  with a linear far field flow.
%The considered examples are related to  the two-phase quadrature domain theory \cite{shah}, \cite{gardiner}. 

Our study showed the possibility for the control  of the interface  via the two-phase mother body for the two-phase Hele-Shaw problem.

\textbf{Acknowledgments.} The authors are grateful to the referee for helpful suggestions.

%%%%%%%%%%%%%%%%%%%%%%%%%%%%%%%%%%
%% \bibliographystyle{amsplain} %%
%% \bibliography{reflection}    %%

\begin{thebibliography}{1}

\bibitem{mus}
Muskat, M., \emph{Two-fluid systems in porous media. The encroachment
 of water into an oil sand},
Physics \textbf{5},
(1934), 250–-264.



\bibitem{Vas2009}
Vasil'ev, A., \emph{From the Hele-Shaw experiment to integrable systems: a historical overview}, Compl. Anal. Oper. Theory \textbf{3}, no. 2,(2009)  551--585.

\bibitem{Vas2015}
Gustafsson, B., Teodorescu, R., and  Vasil'ev, A., \emph{Classical and stochastic Laplacian growth}, Birkhäuser Verlag, (2015), 315 pp. 

%\bibitem{ST} 
%Saffman, P.G. and Taylor, G.I., \emph{The penetration of a fluid into a porous medium or Hele-Shaw cell containing a more viscous liquid}, Proc. R. Soc. Lond. A \textbf{245}, 312--329.




\bibitem{howison2000}
Howison, S.D., \emph{A note on the two-phase Hele-Shaw problem}, J. Fluid Mech., \textbf{409}, 243--249.

\bibitem{JS}
Jacquard, P. and S\'eguier, P., \emph{Mouvement de deux fluides en contact dans un milieu poreux},
J. de Mec. \textbf{1} (1962), 367--394.




\bibitem{FT}
Friedman A. and Tao, Y., \emph{Nonlinear stability of the Muskat problem with capillary pressure at the free boundary} Nonlinear Anal., \textbf{53} (2003), 45--80.

\bibitem{siegel}
Siegel, M., Caflisch, R.E., and Howison, S., \emph{Global existence, singular solutions,
and ill-posedness for the Muskat problem}, Communications on Pure and Applied Mathematics, Vol. LVII, 
(2004) 0001–-0038.

\bibitem{YT}
Ye, J. and Tanveer, S., \emph{Global solutions for a two-phase Hele-Shaw bubble for a near-circle initial shape}, Compl. Var. Elliptic Eq., \textbf{57} N 1, (2012) 23--61.

\bibitem{crowdy2006}
Crowdy, D., \emph{Exact solutions to the unsteady two-phase Hele-Shaw problem}, Q. J.  Mech. Appl. Maths, \textbf{59}, (2006) 475--485.

\bibitem{dmb}
~Savina, T.V. \& ~Nepomnyashchy, A.A.,  \emph{A dynamical mother body in a Hele-Shaw problem}, Physica D \textbf{240}, (2011)  1156--1163. 

\bibitem{external}
~Savina, T.V. \& ~Nepomnyashchy, A.A.,  \emph{The shape control of a growing air bubble in a Hele-Shaw cell}, SIAM J. Appl. Math. \textbf{75}, (2015)  1261--1274. 



\bibitem{HowBubble}
~Howison, S.D.   \emph{Bubble growth in porous media and Hele-Shaw cells}, Proceedings of the Royal Society of Edinburgh, 
\textbf{102A},  (1986)  141--148.

\bibitem{di}
~Di Benedetto, E. and ~Friedman, A.  \emph{Bubble growth in porous media}, Indiana University Mathematics Journal, \textbf{35},  No. 3, (1986) 573--606.

\bibitem{karp1}
~Karp, L.,  \emph{Asymptotic properties of unbounded quadrature domains in the plane}, European Journal of Applied Mathematics, \textbf{26}, No. 3, (2015) 175--191.

\bibitem{karp2}
~Karp, L.,  \emph{Generalized Newtonian potential and its applications}, Journal of  Mathematical Analysis and its Applications, \textbf{174},  (1993) 480--497.


\bibitem{shah}
~Emamizadeh, B., Prajapat, J.V., and ~Shahgholian, H.,  \emph{A two phase free boundary problem related to quadrature domains}, Potential Anal., \textbf{34},  (2011) 119--138.

\bibitem{gardiner}
~Gardiner, S.J and ~Sj\"odin, T.  \emph{Two-phase quadrature domains},
J. D'Analyse Mathematique, \textbf{116}, N 1, (2012) 335--354.







%\bibitem{langer} Langer, J.S., \emph{Instabilities and pattern formation in crystal growth}, Rev. Mod. Phys. \textbf{52}, (1980) 1-28.

%\bibitem{sawada} Sawada Y., Dougherty A., and Gollub J.P., \emph{Dendritic and fractal patterns in electrolytic metal deposits}, Phys. Rev. Lett. \textbf{56} no. 12, (1986) 1260-1263.


%\bibitem{bates}
%~Alikakos, N.,~Bates, P., and ~ Chen, X., \emph{ Convergence of the
%Cahn-Hilliard equation to the Hele-Shaw model}, Arch. Rat. Mech.
%Anal., \textbf{128} (1994), 165--205.

%\bibitem{pego}
%~Pego, R.L., \emph{Front migration in the nonlinear Cahn-Hilliard
%equation}, Proc. Roy. Soc. London A \textbf{422} (1989), 261--278.

%\bibitem{suppr}
%Savina, T.V., Nepomnyashchy, A.A., Brandon, S., Lewin D.R.,and Golovin, A.A., \emph{Suppressing morphological instability via feedback control}, J. Crystal Growth, \textbf{240} (2002), 292-304.

%\bibitem{paterson}
%~Paterson, L.  \emph{Radial fingering in a Hele-Shaw cell}, J. Fluid Mech., 
%\textbf{113},  (1981) 513--529.


%\bibitem{170}
%~Hohlov, Yu.E., ~Howison, S.D., ~Huntingford, C., ~Ockendon, J.R. \& ~Lacey, A.A.  \emph{A model for non-smooth free-boundaries in Hele-Shaw flows}, Quart. J. Mech. Appl. Math., \textbf{47},  (1994) 107--128.


%\bibitem{117}
%~Galin, L.A.,   \emph{Unsteady filtration with a free surface}, Dokl. Acad. Nauk USSR, 
%\textbf{47},  (1945) 246--249 (in Russian).


%\bibitem{262}
%~Polubarinova-Kochina, P.Ya.  \emph{On a problem of the motion of the contour of a petroleum shell}, Prikl. Matem. Mech, 
%\textbf{9},   No. 1, (1945) 79--90 (in Russian).

%\bibitem{263}
%~Polubarinova-Kochina, P.Ya.   \emph{Concerning unsteady motions in the theory of filtration}, Dokl. Acad. Nauk USSR, 
%\textbf{47},   No. 4, (1945) 254--257 (in Russian).


%\bibitem{GusVas}
%~Gustafsson, B. \& ~Vasil'ev, A.  \emph{Conformal and potential analysis in Hele-Shaw cells},
%Birkh\"auser, 2006.




\bibitem{davis}
~Davis, Ph.  \emph{The Schwarz function and its applications},
Carus Mathematical Monographs, MAA, 1979.


\bibitem{laguna}
~Khavinson, D.  \emph{Holomorphic partial differential equations and classical potential theory},
Universidad de La Laguna, 1996.



\bibitem{nonloc}
~Savina, T.  \emph{On non-local reflection for  elliptic equations of the second
order in $\mathbb{R}^2$
(the Dirichlet condition)}, Trans. Amer. Math. Soc. \textbf{364}, no. 5, (2012) 2443-2460.

\bibitem{shapiro}
~Shapiro, H.S.  \emph{The Schwarz function and its generalization to higher dimensions}, John Wiley and Sons, Inc., 1992.


\bibitem{howison92}
~Howison, S.D.  \emph{Complex variable methods in Hele-Shaw moving boundary problems}, European J. Appl. Math., 
\textbf{3},   No. 3, (1992) 209--224.



\bibitem{cummings}
Cummings, L.J., Howison, S.D. \& ~King, J.R.  \emph{Two-dimensional Stokes and Hele-Shaw flows with free surfaces}, J. Appl. Math., 
\textbf{10},  (1999) 635--680.


\bibitem{KhMP}
~Khavinson, D., ~Mineev-Weinstein, M. \& ~Putinar, M.  \emph{Planar elliptic growth}, Complex Analysis and Operator Theory,
\textbf{3}, No. 2, (2009) 425--451.

\bibitem{lacey}
~Lacey, A.A.   \emph{Moving boundary problems in the flow of liquid through porous media}, J. Austral. Math. Soc., 
\textbf{B24}, (1982) 171--193.

\bibitem{Gu1}
~Gustafsson, B.  \emph{On mother bodies of convex polyhedra},
SIAM  J. Math. Anal., \textbf{29}, N 5, (1998) 1106--1117.

\bibitem{GuSa2}
~Gustafsson, B. and ~Sakai, M.  \emph{On potential theoretic skeletons of polyhedra},
Geometriae Dedicata, \textbf{76}, (1999) 1--30.


\bibitem{sss2005}
~Savina, T.V., ~Sternin, B.Yu. \& ~Shatalov, V.E.  \emph{On a minimal element for a family of bodies producing the same external gravitational field}, Appl. Anal. \textbf{84},  no. 7, (2005) 649-668. 




\bibitem{algtop}
Hatcher, A. \emph{Algebraic topology}, Cambridge University Press, 2002.


\bibitem{BE}
Bateman, H. \& Erd\'elyi, A., \emph{Higher transcendental functions}, MC Grow-Hill Book Company, 1955.



%\bibitem{Zid1}
%~Zidarov, D.  \emph{On solution of some inverse problems for potential fields and its application to questions 
%in geophysics},
%Publ. House of Bulg. Acad. of Sci., Sofia,  (Russian), 1968.

%\bibitem{Zid2}
%~Zidarov, D.  \emph{Inverse gravimetric  problem 
%in geoprospecting and geodesy}, Elsevier, Amsterdam, 1990.



%\bibitem{ts}
%~Tsirulskii, A.V.  \emph{Functions of complex variable in the theory and methods of potential geophysical fields}, URO AN SSSR, Sverdlovsk,  (Russian), 1990.



%\bibitem{willis}
%Willis, N.J. \emph{Bistatic  radar}, SciTech, 2005.

\bibitem{ww}
~Whittaker, E.T. and ~Watson, G.N. \emph{A course of modern analysis},  Cambridge University Press, 1996.

%\bibitem{VascKad}
%~Vasconcelos, G.L. and ~Kadanoff, L.P. \emph{Stationary solutions for the Saffman-Taylor problem with surface tension}, Phys. Rev. A, \textbf{44}, (1991),  6490-6460.

%\bibitem{Vasc}
% ~Vasconcelos, G.L. \emph{Exact solutions for Hele-Shaw flows with surface tension: the Schwarz-function approach}, Phys. Rev. E, \textbf{48}, (1993), No. 2, R658-R660.



%\bibitem{GusShap}
%~Gustafsson, B. and ~Shapiro, H.S.  \emph{What is a quadrature domain?}, in P.~Ebenfelt, B.~Gustafsson, D.~Khavinson, M.~Putinar (editors):  \emph{Quadrature domains and applications}, a Harold S. Shapiro Anniversary Volume,
%Birkh\"auser, 2005.



%\bibitem{sakai}
%~Sakai, M.   \emph{Null quadrature domains}, J. Analyse. Math., 
%\textbf{40}, (1981) 144--154.






%\bibitem{81}
%E.~Di Benedetto and A.~Friedman, \emph{The ill-posed Hele-Shaw model and the Stephan problem for supercooled water}, Trans. Amer. Math. Soc., \textbf{282},  (1984), %No. 1, 183--204.\bibitem{cum}










%\bibitem{Gus}
%~Gustafsson, B. (2004) \emph{Lectures on balayage}, in Clifford algebras and potential theory, 17-63,
%Univ. Joensuu, Joensuu.







%\bibitem{177}
%~Howison, S.D.  (1986) \emph{Cusp development in Hele-Shaw flow with a free surface}, SIAM J. Appl. Math., 
%\textbf{46},   No. 1, 20--26.









%\bibitem{mcdonald}
%~McDonald, N.R. (2011) \emph{Generalized Hele-Shaw flow: A Schwarz function approach}, European J. Appl. Math., 
%\textbf{22},   517--532.






%\bibitem{bihar}
%T.V.~Savina,  \emph{On the dependence of the reflection operator on boundary conditions for biharmonic functions}, J. Math. Anal. Appl., \textbf{370} (2010), 716-725. 

%\bibitem{sss95}
%T.V.~Savina, B.Yu.~Sternin and V.E.~Shatalov,  \emph{Notes on mother body in geophysics}, Preprint Max-Plank Institut fur Mathematik, Bonn, MPI/95-90, 23 p. 





%\bibitem{KhMPT}
%D.~Khavinson, M.~Mineev-Weinstein, M.~Putinar, and R.~Teodorescu \emph{Lemniscates do not %survive Laplacian growth}, Mathematical Research Letters, 
%\textbf{17},  (2010), No. 2, 337--343.

\end{thebibliography}
%%%%%%%%%%%%%%%%%%%%%%%%%%%%%%%%%%

\providecommand{\bysame}{\leavevmode\hbox to3em{\hrulefill}\thinspace}

\end{document}